\documentclass[leqno]{amsart}

\usepackage[leqno]{amsmath}
\usepackage{amssymb,amsthm}
\usepackage{enumerate, latexsym, mathrsfs, mdwlist}

\usepackage[final,factor=1100,stretch=10,shrink=10]{microtype}

\SetProtrusion{encoding={*},family={bch},series={*},size={6,7}}
              {1={ ,750},2={ ,500},3={ ,500},4={ ,500},5={ ,500},
               6={ ,500},7={ ,600},8={ ,500},9={ ,500},0={ ,500}}

\SetExtraKerning[unit=space]
    {encoding={*}, family={bch}, series={*}, size={footnotesize,small,normalsize}}
    {\textendash={400,400}, 
     "28={ ,150}, 
     "29={150, }, 
     \textquotedblleft={ ,150}, 
     \textquotedblright={150, }} 

\swapnumbers
\theoremstyle{definition}
\newtheorem{nul}{}[section]
\newtheorem{dfn}[nul]{Definition}

\newtheorem{cnstr}[nul]{Construction}
\newtheorem{rem}[nul]{Remark}
\newtheorem{ntn}[nul]{Notation}
\newtheorem{exm}[nul]{Example}

\newtheorem*{dfn*}{Definition}
\newtheorem*{axm*}{Axiom}
\newtheorem*{ntn*}{Notation}
\newtheorem*{exm*}{Example}
\newtheorem*{exr*}{Exercise}
\newtheorem*{int*}{Intuition}
\newtheorem*{qst*}{Question}

\theoremstyle{plain}

\newtheorem{thm}[nul]{Theorem}
\newtheorem{prp}[nul]{Proposition}
\newtheorem{lem}[nul]{Lemma}

\newtheorem{cor}{Corollary}[nul]
\newtheorem*{thm*}{Theorem}
\newtheorem*{prp*}{Proposition}
\newtheorem*{cor*}{Corollary}
\newtheorem*{lem*}{Lemma}
\newtheorem*{cnj*}{Conjecture}

\numberwithin{equation}{nul}

\DeclareMathOperator{\End}{End}

\DeclareMathOperator{\Exc}{Exc}
\DeclareMathOperator{\Ext}{Ext}

\DeclareMathOperator{\Fun}{Fun}
\DeclareMathOperator{\id}{id}

\DeclareMathOperator{\Map}{Map}

\DeclareMathOperator{\Thy}{Thy}

\newcommand{\DD}{\mathbf{D}}

\newcommand{\FF}{\mathbf{F}}

\newcommand{\KK}{\mathbf{K}}

\newcommand{\NN}{\mathbf{N}}
\newcommand{\OO}{\mathbf{O}}

\renewcommand{\SS}{\mathbf{S}}

\newcommand{\D}{\mathrm{D}}

\newcommand{\Add}{\mathrm{Add}}
\newcommand{\Alg}{\mathbf{Alg}}
\newcommand{\Cat}{\mathbf{Cat}}
\newcommand{\Fin}{\mathbf{Fin}}
\newcommand{\Kan}{\mathbf{Kan}}

\newcommand{\Pair}{\mathbf{Pair}}

\newcommand{\Sp}{\mathbf{Sp}}
\newcommand{\VaddWald}{\mathrm{D}_{\mathrm{fiss}}(\mathbf{Wald}_{\infty})}
\newcommand{\VWald}{\mathrm{D}(\mathbf{Wald}_{\infty})}
\newcommand{\Wald}{\mathbf{Wald}}

\newcommand{\fiss}{\mathrm{fiss}}
\newcommand{\op}{\mathrm{op}}

\newcommand{\coloneq}{\mathrel{\mathop:}=}

\makeatletter 
\def\revddots{\mathinner{\mkern1mu\raise\p@ 
\vbox{\kern7\p@\hbox{.}}\mkern2mu 
\raise4\p@\hbox{.}\mkern2mu\raise7\p@\hbox{.}\mkern1mu}} 
\makeatother

\usepackage{tikz}
\usetikzlibrary{matrix,arrows,decorations}

\pgfdeclaredecoration{single line}{initial}{\state{initial}[width=\pgfdecoratedpathlength-1sp]{\pgfpathmoveto{\pgfpointorigin}}\state{final}{\pgfpathlineto{\pgfpointorigin}}} 

\newcommand{\fromto}[2]{{#1}\ \tikz[baseline]\draw[>=stealth,->](0,0.5ex)--(0.5,0.5ex);\ {#2}}
\newcommand{\into}[2]{{#1}\ \tikz[baseline]\draw[>=stealth,right hook->](0,0.5ex)--(0.5,0.5ex);\ {#2}}

\newcommand{\cofto}[2]{{#1}\ \tikz[baseline]\draw[>=stealth,>->](0,0.5ex)--(0.5,0.5ex);\ {#2}}

\newcommand{\equivto}[2]{{#1}\ \tikz[baseline]\draw[>=stealth,->,font=\scriptsize,inner sep=0.5pt](0,0.5ex)--node[above]{$\sim$}(0.5,0.5ex);\ {#2}}
\newcommand{\goesto}[2]{{#1}\ \tikz[baseline]\draw[|->](0,0.5ex)--(0.5,0.5ex);\ {#2}}

\renewcommand{\to}{\ \tikz[baseline]\draw[>=stealth,->](0,0.5ex)--(0.5,0.5ex);\ }

\title{Multiplicative structures on algebraic $K$-theory}

\author{Clark Barwick}
\address{Massachusetts Institute of Technology, Department of Mathematics, Building 2, 77 Massachusetts Avenue, Cambridge, MA 02139-4307, USA}
\email{clarkbar@gmail.com}

\begin{document}

\maketitle


\setcounter{section}{-1}
\section{Introduction} Dan Kan playfully described the theory of $\infty$-categories as \emph{the homotopy theory of homotopy theories}. The aim of this paper, which is a sequel to \cite{K1}, is to show that algebraic $K$-theory is the \emph{stable} homotopy theory of homotopy theories, and it interacts with algebraic structures accordingly. To explain this assertion, let's recap the contents of \cite{K1}.

\begin{nul}\label{nul:Waldinftycat} The kinds of homotopy theories under consideration in this paper are \emph{Waldhausen $\infty$-categories} \cite[Df. 2.7]{K1}. (We employ the quasicategory model of $\infty$-categories for technical convenience.) These are $\infty$-categories with a zero object and a distinguished class of morphisms (called \emph{cofibrations} or \emph{ingressive morphisms}) that satisfies the following conditions.
\begin{enumerate}[(\ref{nul:Waldinftycat}.1)]
\item Any equivalence is ingressive.
\item Any morphism from the zero object is ingressive.
\item Any composite of ingressive morphisms is ingressive.
\item The (homotopy) pushout of an ingressive morphism along any morphism exists and is ingressive.
\end{enumerate}
A pushout of a cofibration $\cofto{X}{Y}$ along the map to the zero object is to be viewed as a \emph{cofiber sequence}
\begin{equation*}
X\ \tikz[baseline]\draw[>=stealth,>->](0,0.5ex)--(0.5,0.5ex);\ Y\ \tikz[baseline]\draw[>=stealth,->](0,0.5ex)--(0.5,0.5ex);\ Y/X.
\end{equation*}

Examples of this structure abound: pointed $\infty$-categories with all finite colimits, exact categories in the sense of Quillen, and many categories with cofibrations and weak equivalences in the sense of Waldhausen all provide examples of Waldhausen $\infty$-categories.

Write $\Wald_{\infty}$ for the $\infty$-category whose objects are Waldhausen $\infty$-categories and whose morphisms are functors that are \emph{exact} in the sense that they preserve the cofiber sequences. This is a compactly generated $\infty$-category \cite[Pr. 4.7]{K1} that admits direct sums \cite[Pr. 4.6]{K1}.
\end{nul}

\begin{nul}\label{nul:iandr} We will be interested in invariants that \emph{split cofiber sequences} in Waldhausen $\infty$-categories. To make this precise, we construct, for any Waldhausen $\infty$-category $C$, let $E(C)$ denote the $\infty$-category of cofiber sequences
\begin{equation*}
X\ \tikz[baseline]\draw[>=stealth,>->](0,0.5ex)--(0.5,0.5ex);\ Y\ \tikz[baseline]\draw[>=stealth,->](0,0.5ex)--(0.5,0.5ex);\ Y/X.
\end{equation*}
This is a Waldhausen $\infty$-category \cite[Pr. 5.11]{K1} in which a morphism
\begin{equation*}
\begin{tikzpicture} 
\matrix(m)[matrix of math nodes, 
row sep=4ex, column sep=4ex, 
text height=1.5ex, text depth=0.25ex] 
{U&V&V/U\\ 
X&Y&Y/X\\}; 
\path[>=stealth,->,font=\scriptsize] 
(m-1-1) edge (m-1-2) 
edge (m-2-1) 
(m-1-2) edge (m-1-3)
edge (m-2-2)
(m-1-3) edge (m-2-3)
(m-2-1) edge (m-2-2)
(m-2-2) edge (m-2-3); 
\end{tikzpicture}
\end{equation*}
is ingressive just in case each of
\begin{equation*}
\cofto{U}{X}\textrm{,\qquad}\cofto{V}{Y}\textrm{,\quad and\quad}\cofto{V/U}{Y/X}
\end{equation*}
is ingressive. We have an exact functor $m\colon\fromto{E(C)}{C}$ defined by the assignment
\begin{equation*}
\goesto{\left[X\ \tikz[baseline]\draw[>=stealth,>->](0,0.5ex)--(0.5,0.5ex);\ Y\ \tikz[baseline]\draw[>=stealth,->](0,0.5ex)--(0.5,0.5ex);\ Y/X\right]}{Y}.
\end{equation*}
We also have an exact functor $i\colon\fromto{C\oplus C}{E(C)}$ defined by the assignment
\begin{equation*}
\goesto{(X,Z)}{\left[X\ \tikz[baseline]\draw[>=stealth,>->](0,0.5ex)--(0.5,0.5ex);\ X\vee Z\ \tikz[baseline]\draw[>=stealth,->](0,0.5ex)--(0.5,0.5ex);\ Z\right]}
\end{equation*}
as well as a retraction $r\colon\fromto{E(C)}{C\oplus C}$ defined by the assignment
\begin{equation*}
\goesto{\left[X\ \tikz[baseline]\draw[>=stealth,>->](0,0.5ex)--(0.5,0.5ex);\ Y\ \tikz[baseline]\draw[>=stealth,->](0,0.5ex)--(0.5,0.5ex);\ Y/X\right]}{(X,Y/X)}.
\end{equation*}
\end{nul}

\begin{nul}\label{nul:add} Now let $\SS$ denote the $\infty$-category of spaces. In \cite[\S 6]{K1}, I constructed a homotopy theory $\D_{\fiss}(\Wald_{\infty})$ such that homology theories (i.e., reduced, $1$-excisive functors)
\begin{equation*}
\fromto{\D_{\fiss}(\Wald_{\infty})}{\SS}
\end{equation*}
are essentially the same data \cite[Th. 7.4]{K1} as functors $\phi\colon\fromto{\Wald_{\infty}}{\SS}$ with the following properties.
\begin{enumerate}[(\ref{nul:add}.1)]
\item $\phi$ is \emph{finitary} in the sense that it preserves filtered colimits.
\item $\phi$ is \emph{reduced} in the sense that $\phi(0)=\ast$.
\item $\phi$ \emph{splits cofiber sequences} in the sense that the exact functor $r$ induces an equivalence
\begin{equation*}
\equivto{\phi(E(C))}{\phi(C)\times\phi(C)}.
\end{equation*}
\item $\phi$ is \emph{grouplike} in the sense that the multiplication
\begin{equation*}
\phi(m\circ i)\colon\fromto{\phi(C)\times\phi(C)\simeq\phi(E(C))}{\phi(C)}
\end{equation*}
defines a grouplike $H$-space structure on $\phi(C)$.
\end{enumerate}
$\D_{\fiss}(\Wald_{\infty})$ is called the \emph{fissile derived $\infty$-category} of Waldhausen $\infty$-categories, and its objects are called \emph{fissile virtual Waldhausen $\infty$-categories}. It is possible to be quite explicit about these objects: fissile virtual Waldhausen $\infty$-categories are functors $\fromto{\Wald_{\infty}^{\op}}{\SS}$ that satisfy the dual conditions to (\ref{nul:add}.1--3). The homotopy theory $\D_{\fiss}(\Wald_{\infty})$ is compactly generated and it admits all direct sums; furthermore, suspension in $\D_{\fiss}(\Wald_{\infty})$ is given by the geometric realization of Waldhausen's $S_{\bullet}$ construction \cite[Cor. 6.9.1]{K1}.
\end{nul}

\begin{nul}\label{nul:aboutP1} From any finitary reduced functor $F\colon\fromto{\D_{\fiss}(\Wald_{\infty})}{\SS}$ one may extract the Goodwillie differential $P_1F$, which is the nearest excisive approximation to $F$, or, in other words, the best approximation to $F$ by a homology theory \cite{MR2026544}. This approximation is given explicitly by the colimit of the sequence
\begin{equation*}
F\to\Omega\circ F\circ\Sigma\to\cdots\to\Omega^n\circ F\circ\Sigma^n\to\cdots.
\end{equation*}

Since $P_1F$ is excisive, it factors naturally through the functor $\Omega^{\infty}\colon\fromto{\Sp}{\SS}$. Consequently, for any fissile virtual Waldhausen $\infty$-category $X$, we obtain a homology theory
\begin{equation*}
P_1F(X)\colon\fromto{\SS_{\ast}}{\SS}.
\end{equation*}
Unwinding the definitions, we find that this homology theory is itself the Goodwillie differential of the functor
\begin{equation*}
\goesto{T}{F(T\otimes X)},
\end{equation*}
where $\otimes$ denotes the tensor product
\begin{equation*}
\fromto{\SS_{\ast}\times\D_{\fiss}(\Wald_{\infty})}{\D_{\fiss}(\Wald_{\infty})}
\end{equation*}
guaranteed by the identification of presentable pointed $\infty$-categories with modules over $\SS_{\ast}$ \cite[Pr. 6.3.2.11]{HA}.
\end{nul}

\begin{nul} The main result of \cite[\S 10]{K1} can now be stated as follows. We have a Waldhausen $\infty$-category $\Fin_{\ast}$ of pointed finite sets, in which the cofibrations are injective (pointed) maps; if $I\colon\fromto{\D_{\fiss}(\Wald_{\infty})}{\SS}$ denotes evaluation at $\Fin_{\ast}$ (so that $I(X)=X(\Fin_{\ast})$), then algebraic $K$-theory may be identified as
\begin{equation*}
K\simeq P_1I.
\end{equation*}
This gives a ``local'' universal property for algebraic $K$-theory: for any fissile virtual Waldhausen $\infty$-category $X$, the homology theory $\KK(X)$ is the Goodwillie differential of the functor
\begin{equation*}
\goesto{T}{I(T\otimes X)}.
\end{equation*}
Our version of Waldhausen's Additivity Theorem states that this differential only requires a single delooping: $\KK(X)$ is simply the homology theory $\fromto{\SS_{\ast}}{\SS}$ given by
\begin{equation*}
\goesto{T}{\Omega I(\Sigma T\otimes X)},
\end{equation*}
or, equivalently, since suspension in $\D_{\fiss}(\Wald_{\infty})$ is given by Waldhausen's $S_{\bullet}$ construction, the assignment
\begin{equation*}
\goesto{T}{\Omega I(T\otimes S_{\bullet}(X))},
\end{equation*}
\end{nul}

\begin{nul} In this paper, I go further. I construct (Pr. \ref{prp:Waldinftyotimes}) a symmetric monoidal structure on $\Wald_{\infty}$ in which the tensor product $C\otimes D$ represents ``bi-exact'' functors
\begin{equation*}
\fromto{C\times D}{E},
\end{equation*}
i.e., functors that preserve cofiber sequences separately in each variable. The unit therein is simply the Waldhausen $\infty$-category $\Fin_{\ast}$. I then descend this symmetric monoidal structure to one on $\D_{\fiss}(\Wald_{\infty})$ with the property that it preserves colimits separately in each variable (Pr. \ref{prp:Dfissotimes}).

Now the functor represented by $\Fin_{\ast}$ is the unit for the Day convolution symmetric monoidal structure constructed by Saul Glasman \cite{GlasmanDay} on the $\infty$-category of functors from $\DD_{\fiss}$ to spaces. Its differential --- i.e., algebraic $K$-theory --- is therefore the unit among homology theories on $\D_{\fiss}(\Wald_{\infty})$. That is, it plays precisely the same role among homology theories on $\D_{\fiss}(\Wald_{\infty})$ that is played by the sphere spectrum in the $\infty$-category of spectra. It therefore has earned the mantle \emph{the stable homotopy theory of homotopy theories}.

Just as one may describe the stable homotopy groups of a pointed space $X$ as Ext groups out of the unit:
\begin{equation*}
\pi_n^{s}(X)\cong\Ext^{-n}(\Sigma^{\infty}S^0,\Sigma^{\infty}X)
\end{equation*}
in the stable homotopy category, so too may one describe the algebraic $K$-theory groups of a Waldhausen $\infty$-category $C$ as Ext groups out of the unit:
\begin{equation*}
K_n(C)\cong\Ext^{-n}(\Sigma^{\infty}\Fin_{\ast},\Sigma^{\infty}C)
\end{equation*}
in the \emph{stable homotopy category of Waldhausen $\infty$-categories}.
\end{nul}

\begin{nul} Algebraic $K$-theory is therefore naturally multiplicative, and so it inherits homotopy-coherent algebraic structures on Waldhausen $\infty$-categories. That is, I show (Cor. \ref{cor:Kpreservesalg}) that if an $\infty$-operad $O$ acts on a Waldhausen $\infty$-category $C$ via functors that are exact separately in each variable, then there is an induced action of $O$ on both the space $K(C)$ and the spectrum $\KK(C)$. As a corollary (Ex. \ref{ex:Deligne}), I deduce that for any $1\leq n\leq\infty$, the algebraic $K$-theory of an $E_n$-algebra in a suitable symmetric monoidal $\infty$-category is an $E_{n-1}$ ring spectrum. In particular, we note that the $K$-theory of an $E_n$ ring is an $E_{n-1}$ ring, and the $A$-theory of any $n$-fold loopspace is an $E_{n-1}$ ring spectrum. These sorts of results are $K$-theoretic analogues of the so-called (homological) \emph{Deligne Conjecture} \cite{MR1890736,MR1805894}.
\end{nul}

\begin{nul} In fact, the main result (Th. \ref{thm:main}) is much more general: I actually show that any additive theory that can be expressed as the additivization of a \emph{multiplicative} (i.e., multiplicative) theory is itself multiplicative in a canonical fashion. This yields a uniform way of reproducing conjectures ``of Deligne type'' for theories that are both additive and multiplicative.

Of course since algebraic $K$-theory is the stable homotopy theory of Waldhausen $\infty$-categories, it is \emph{initial} as a theory that is both additive and multiplicative. The sphere spectrum is the initial $E_{\infty}$ object of the stable $\infty$-category of spectra; similarly, the object of the stable $\infty$-category of Waldhausen $\infty$-categories that represents algebraic $K$-theory is the initial $E_{\infty}$ object. Just as the universal property of algebraic $K$-theory \cite{K1} gives a uniform construction of trace maps, this result gives a uniform construction of \emph{multiplicative} trace maps --- in particular, any additive and multiplicative theory accepts a \emph{unique} (up to a contractible choice) multiplicative trace map (Th. \ref{multtrace}).

The passage to higher categories is a \emph{sine qua non} of this result. Indeed, note that when $n\geq3$, an $E_n$ structure on an ordinary category is tantamount to a symmetric monoidal structure. As a result, it is difficult to identify, e.g., $E_3$ structures on the $K$-theory spectrum of an $E_4$-algebra without employing some higher categorical machinery.
\end{nul}

\begin{rem} Note that the form of algebraic $K$-theory we study has an exceptionally strong compatibility with the tensor product. For example, the Barratt--Priddy--Quillen--Segal theorem implies the endomorphism spectrum $\End(\Sigma^{\infty}\Fin_{\ast})$ is the sphere spectrum. That is, the form of algebraic $K$-theory studied here is strongly unital. This is of course \emph{false} for any form of algebraic $K$-theory that applies only to Waldhausen $\infty$-categories in which every morphism is ingressive. I intend to return to this observation in future work.
\end{rem}

\begin{rem} The heart of the proof is to use the description of $K$-theory as a Goodwillie differential. The result actually follows from a quite general fact about the interaction between the Goodwillie calculus and symmetric monoidal structures. Namely, the Goodwillie differential of a multiplicative functor between suitable symmetric monoidal $\infty$-categories inherits a canonical multiplicative structure. This fact, which may be of independent interest, doesn't seem to be recorded anywhere in the literature. So I do so in this paper (Pr. \ref{prp:diffofalg}).
\end{rem}

\begin{rem} Some variants of some of these results can be found in the literature.

Elmendorf and Mandell \cite{MR2254311} constructed an algebraic $K$-theory of multicategories, which the show is lax monoidal as a functor to symmetric spectra. Consequently, they deduce that any operad that acts on a permutative category will also act on its $K$-theory.

In later work of Blumberg and Mandell \cite[Th. 2.6]{MR2805994}, it is shown that the $K$-theory of an ordinary Waldhausen category equipped with an action of a categorical operad inherits an action of the nerve of that operad.

In independent work, Blumberg, Gepner, and Tabuada \cite{BGT2} have proved that algebraic $K$-theory is initial among additive and multiplicative functors from an $\infty$-category of idempotent complete stable $\infty$-categories. (They used this to uniquely characterize the cyclotomic trace.) The result here shows that algebraic $K$-theory is initial among additive and multiplicative functors on \emph{all} Waldhausen $\infty$-categories.
\end{rem}

\begin{rem} Finally, I emphasize that work of Saul Glasman \cite{GlasmanDay} made it possible for me to sharpen the results of this paper significantly. Previous versions of this paper did not contain the full strength of the universality of algebraic $K$-theory as an additive and multiplicative theory.
\end{rem}


\section{Tensor products of Waldhausen $\infty$-categories} The first thing we need to understand is the symmetric monoidal structure on the $\infty$-category of Waldhausen $\infty$-categories. As in \cite[\S 4]{K1}, we will regard $\Wald_{\infty}$ as formally analogous to the nerve of the ordinary category $V(k)$ of vector spaces over a field $k$. The tensor product $V\otimes W$ of vector spaces is defined as the vector space that represents \emph{multilinear} maps $\fromto{V\times W}{X}$, i.e., maps that are linear separately in each variable. In perfect analogy with this, the tensor product $ C\otimes D$ of Waldhausen $\infty$-categories is defined as the Waldhausen $\infty$-category that represents functors $\fromto{ C\times D}{ E}$ that are exact separately in each variable.

\begin{ntn}  Let $\Lambda(\FF)$ denote the following ordinary category. The objects will be finite sets, and a morphism $\fromto{J}{I}$ will be a map $\fromto{J}{I_{+}}$; one composes $\psi\colon\fromto{K}{J_{+}}$ with $\phi\colon\fromto{J}{I_{+}}$ by forming the composite
\begin{equation*}
K\ \tikz[baseline]\draw[>=stealth,->,font=\scriptsize](0,0.5ex)--node[above]{$\psi$}(0.5,0.5ex);\ J_{+}\ \tikz[baseline]\draw[>=stealth,->,font=\scriptsize](0,0.5ex)--node[above]{$\phi_{+}$}(0.5,0.5ex);\ I_{++}\ \tikz[baseline]\draw[>=stealth,->,font=\scriptsize](0,0.5ex)--node[above]{$\mu$}(0.5,0.5ex);\ I_{+},
\end{equation*}
where $\mu\colon\fromto{I_{++}}{I_{+}}$ is the map that simply identifies the two added points. (Of course $\Lambda(\FF)$ is equivalent to the category $\mathbf{Fin}_{\ast}$ of pointed finite sets, but we prefer to think of the objects of $\Lambda(\FF)$ as unpointed. This is the natural perspective on this category from the theory of operator categories \cite{opcat}.)

For any morphism $\phi:\fromto{J}{I}$ of $\Lambda(\FF)$ and any $i\in I$, write $J_i$ for the fiber $\phi^{-1}(\{i\})$.
\end{ntn}

\begin{nul} One way to write down a symmetric monoidal $\infty$-category \cite[Ch. 2]{HA} is to give the data of the space of maps out of any tensor product of any finite collection of objects. More precisely, a \emph{symmetric monoidal $\infty$-category} is a cocartesian fibration $p:\fromto{C^{\otimes}}{N\Lambda(\FF)}$ such that for any finite set $I$, the various maps $\chi_i:\fromto{I}{\{i\}_+}$ such that $\chi_i^{-1}(\{i\})=\{i\}$ together specify an equivalence of $\infty$-categories
\begin{equation*}
\equivto{C^{\otimes}_I}{\prod_{i\in I}C^{\otimes}_{\{i\}}}.
\end{equation*}
The objects of $C^{\otimes}$ are, in effect, $(I,X_I)$ consisting of a finite set $I$ and a collection $X_I=\{X_i\}_{i\in I}$ of objects of $C$. Morphisms $\fromto{(J,Y_J)}{(I,X_I)}$ of $C^{\otimes}$ are essentially pairs $(\omega,\phi_I)$ consisting of morphisms $\omega:\fromto{J}{I}$ of $\Lambda(\FF)$ and families of morphisms
\begin{equation*}
\left\{\phi_i:\fromto{\bigotimes_{j\in J_i}Y_j}{X_i}\right\}_{i\in I}.
\end{equation*}
\end{nul}

\begin{exm} For any $\infty$-category that admits all finite products, there is a corresponding symmetric monoidal $\infty$-category $C^{\times}$ called the \emph{cartesian} symmetric monoidal $\infty$-category.
\end{exm}

We will be particularly interested in identifying a suitable subcategory of the $\infty$-category $\Pair_\infty^{\times}$, where $\Pair_{\infty}$ denotes the $\infty$-category of pairs of $\infty$-categories \cite[Df. 1.11]{K1}.

\begin{dfn} Suppose $I$ a finite set, and suppose $ C_I\coloneq( C_i)_{i\in I}$ an $I$-tuple of Waldhausen $\infty$-categories. For any Waldhausen $\infty$-category $ D$, a functor of pairs $\fromto{\prod C_I}{ D}$ is said to be \textbf{\emph{exact separately in each variable}} if, for any element $i\in I$ and any collection of objects $(X_j)_{j\in I-\{i\}}\in\prod C_{I-\{i\}}$, the functor
\begin{equation*}
 C_i\cong C_i\times\prod_{j\in I-\{i\}}\{X_j\}\ \tikz[baseline]\draw[>=stealth,right hook->](0,0.5ex)--(0.5,0.5ex);\ \prod C_I\ \tikz[baseline]\draw[>=stealth,->](0,0.5ex)--(0.5,0.5ex);\  D
\end{equation*}
carries cofibrations to cofibrations and is exact as a functor of pairs between Waldhausen $\infty$-categories.
\end{dfn}

\begin{nul} Note that we do not assume the \emph{cubical cofibrancy criterion} that appears in Blumberg--Mandell \cite[Df. 2.4]{MR2805994}. It seems that the authors of this paper required it to guarantee a compatibility with the ``all at once'' iterated $S_{\bullet}$ construction. We will not use such a construction here; the only compatibility we will need to find with the $S_{\bullet}$ construction is Pr. \ref{Lfissismultiplicative}, which deals with one tensor factor at a time. 
\end{nul}

\begin{ntn} Denote by $\Wald_{\infty}^{\otimes}\subset\Pair_{\infty}^{\times}$ the following subcategory. The objects of $\Wald_{\infty}^{\otimes}$ are those objects $(I, C_I)$, where for any $i\in I$, the pair $ C_i$ is a Waldhausen $\infty$-category. A morphism $\fromto{(J, D_J)}{(I, C_I)}$ of $\Pair_{\infty}^{\times}$ is a morphism of $\Wald_{\infty}^{\otimes}$ if and only if, for every element $i\in I$, the functor
\begin{equation*}
\fromto{\prod D_{J_i}}{ C_i}
\end{equation*}
is exact separately in each variable.
\end{ntn}

We now identify the tensor product of Waldhausen $\infty$-categories.

\begin{lem}\label{lem:tensorprodpreWald} Suppose $I$ a finite set, and suppose $ C_I\coloneq( C_i)_{i\in I}$ an $I$-tuple of Waldhausen $\infty$-categories. Then there exist a Waldhausen $\infty$-category $\bigotimes C_I$ and a functor of pairs
\begin{equation*}
f\colon\fromto{\prod C_I}{\bigotimes C_I}
\end{equation*}
such that for every Waldhausen $\infty$-category $ D$, composition with $f$ induces an equivalence between the $\infty$-category $\Fun_{\Wald_{\infty}}(\bigotimes C_I, D)$ and the full subcategory of $\Fun_{\Pair_{\infty}}(\prod C_I, D)$ spanned by the functors of pairs that are exact separately in each variable.
\begin{proof} We construct $\bigotimes C_I$ as a colimit in $\Wald_{\infty}$ in the following manner. First, recall that the forgetful functor $\fromto{\Wald_{\infty}}{\Cat_{\infty}}$ admits a left adjoint $W$. Consider the pushout $K=(I\times\Delta^{1})\cup^{I\times\Delta^{\{0\}}}(I\times\Delta^{\{0\}})^{\rhd}$ and the obvious functor
\begin{equation*}
F\colon\fromto{K}{\Wald_{\infty}}
\end{equation*}
that carries each object of the form $(i,0)$ to the coproduct
\begin{equation*}
\coprod_{(X_j)_{j\in I-\{i\}}\in\prod_{j\in I-\{i\}} C_j}W( C_i),
\end{equation*}
each object of the form $(i,1)$ to the coproduct
\begin{equation*}
\coprod_{(X_j)_{j\in I-\{i\}}\in\prod_{j\in I-\{i\}} C_j} C_i,
\end{equation*}
and the cone point $+\infty$ to $W(\prod_{i\in I} C_i)$. Now the desired Waldhausen $\infty$-category $\bigotimes C_I$ can be constructed as the colimit of $F$.
\end{proof}
\end{lem}

\begin{nul} The construction of this proof is of course the natural analogue of the construction of tensor products of abelian groups. From this description, it is clear that when $I=\varnothing$, then the Waldhausen $\infty$-category $\otimes^0\simeq W(\Delta^0)\simeq N\mathbf{Fin}_{\ast}$, the nerve of the ordinary category of finite pointed sets, in which the cofibrations are the monomorphisms.
\end{nul}

\begin{prp}\label{prp:Waldinftyotimes} The functor $\fromto{\Wald_{\infty}^{\otimes}}{N\Lambda(\FF)}$ is a symmetric monoidal $\infty$-category.
\begin{proof} We first claim that the functor $p\colon\fromto{\Wald_{\infty}^{\otimes}}{N\Lambda(\FF)}$ is a cocartesian fibration; it is an inner fibration because $\Wald_{\infty}^{\otimes}$ is an $\infty$-category \cite[Pr. 2.3.1.5]{HTT}. Now suppose $\phi\colon\fromto{J}{I}$ an edge of $N\!\Lambda$, and suppose $ D_J$ a $J$-tuple of pairs of $\infty$-categories. We want to find a $p$-cocartesian edge of $\Pair_{\infty}^{\otimes}$ covering $\phi$. For this, for any $i\in I$, consider a pair $\bigotimes D_{J_i}$ along with a functor $\fromto{\prod D_{J_i}}{\bigotimes D_{J_i}}$ satisfying the property described in Lemma \ref{lem:tensorprodpreWald}. These fit together to yield a morphism
\begin{equation*}
\fromto{(J, D_J)}{\left(I,\Big(\bigotimes D_{J_i}\Big)_{i\in I}\right)}
\end{equation*}
of $\Wald_{\infty}^{\otimes}$ covering $\phi$. The property described in Lemma \ref{lem:tensorprodpreWald} guarantees that this a locally $p$-cocartesian edge of $\Wald_{\infty}^{\otimes}$, so $p$ is a locally cocartesian fibration. Now to conclude that $p$ is a cocartesian fibration, it is enough to note that for any $2$-simplex
\begin{equation*}
\begin{tikzpicture} 
\matrix(m)[matrix of math nodes, 
row sep=4ex, column sep=4ex, 
text height=1.5ex, text depth=0.25ex] 
{(K, E_K)&&(I, C_I),\\ 
&(J, D_J)&\\}; 
\path[>=stealth,->,font=\scriptsize] 
(m-1-1) edge node[above]{$h$} (m-1-3) 
edge node[below left]{$g$} (m-2-2) 
(m-2-2) edge node[below right]{$f$} (m-1-3); 
\end{tikzpicture}
\end{equation*}
if $f$ and $g$ are locally $p$-cocartesian edges, then so is $h$; this follows directly from our construction of $\bigotimes C_I$.
\end{proof}
\end{prp}

\begin{prp}\label{prp:tensorofWaldgoodforsums} The tensor product functor
\begin{equation*}
\otimes\colon\fromto{\Wald_{\infty}\times\Wald_{\infty}}{\Wald_{\infty}}
\end{equation*}
preserves filtered colimits and direct sums separately in each variable.
\begin{proof} Suppose $ C$ a Waldhausen $\infty$-category. We will show that the endofunctor $-\otimes C$ preserves filtered colimits and direct sums.

Suppose $\Lambda$ a filtered simplicial set, and suppose $ D\colon\fromto{\Lambda^{\rhd}}{\Wald_{\infty}}$ a colimit diagram. Since filtered colimits in $\Wald_{\infty}$ are preserved under the forgetful functor $\fromto{\Wald_{\infty}}{\Pair_{\infty}}$, it follows that for any Waldhausen $\infty$-category $ A$, there is an equivalence of $\infty$-categories
\begin{equation*}
\equivto{\Fun_{\Pair_{\infty}}\left( D_{+\infty}\times C, A\right)}{\lim_{\alpha\in\Lambda}\Fun_{\Pair_{\infty}}\left( D_{\alpha}\times C, A\right)}.
\end{equation*}
To show that the tensor product preserves filtered colimits separately in each variable, it remains to note that, under this equivalence, a functor $\fromto{ D_{+\infty}\times C}{ A}$ that is exact separately in each variable correspond to compatible families of functors $\fromto{ D_{\alpha}\times C}{ A}$ that are exact separately in each variable.

Note that if $0$ is the zero Waldhausen $\infty$-category, then for any Waldhausen $\infty$-category $ A$, any functor $\fromto{0\times C}{ A}$ that is exact separately in each variable is essentially constant, whence $0\otimes C\simeq 0$. Moreover, if $ D$ and $ D'$ are Waldhausen $\infty$-categories, then the two inclusions
\begin{equation*}
\into{ D\times C}{( D\oplus D')\times C}\textrm{\quad and\quad}\into{ D'\times C}{( D\oplus D')\times C}
\end{equation*}
given by $\goesto{(y,x)}{(y,0,x)}$ and $\goesto{(y',x)}{(0,y',x)}$ together induce a functor
\begin{equation*}
\begin{tikzpicture} 
\matrix(m)[matrix of math nodes, 
row sep=3ex, column sep=3ex, 
text height=1.5ex, text depth=0.25ex] 
{\Fun_{\Pair_{\infty}}(( D\oplus D')\times C, A)\\ 
\Fun_{\Pair_{\infty}}( D\times C, A)\times\Fun_{\Pair_{\infty}}( D'\times C, A).\\}; 
\path[>=stealth,->,font=\scriptsize] 
(m-1-1) edge (m-2-1); 
\end{tikzpicture}
\end{equation*}
On the other hand, the coproduct induces a functor
\begin{equation*}
\begin{tikzpicture} 
\matrix(m)[matrix of math nodes, 
row sep=3ex, column sep=3ex, 
text height=1.5ex, text depth=0.25ex] 
{\Fun_{\Pair_{\infty}}( D\times C, A)\times\Fun_{\Pair_{\infty}}( D'\times C, A)\\ 
\Fun_{\Pair_{\infty}}(( D\oplus D')\times C, A).\\}; 
\path[>=stealth,->,font=\scriptsize] 
(m-1-1) edge (m-2-1); 
\end{tikzpicture}
\end{equation*}
It is not hard to see that these functors carry (pairs of) functors that are exact in each variable separately to (pairs of) functors that are exact in each variable separately. Moreover, if $F\colon\fromto{( D\oplus D')\times C}{ A}$ is exact separately in each variable, then
\begin{equation*}
F\simeq F|_{ D\times C}\vee F|_{ D'\times C},
\end{equation*}
and if $G\colon\fromto{ D\times C}{ A}$ and $G'\colon\fromto{ D'\times C}{ A}$ are each exact separately in each variable, then 
\begin{equation*}
G\simeq (G\vee G')|_{ D\times C}\textrm{\quad and\quad}G'\simeq (G\vee G')|_{ D'\times C}.
\end{equation*}
Hence these two functors exhibit an equivalence
\begin{equation*}
( D\oplus D')\otimes C\simeq( D\otimes C)\oplus( D'\otimes C),
\end{equation*}
as desired.
\end{proof}
\end{prp}

\begin{nul} For any integer $m\geq 0$ and any Waldhausen $\infty$-category $ C$, write $m C$ for the iterated direct sum $ C\oplus\cdots\oplus C$, and write $ C^m$ for the iterated tensor product $ C\otimes\cdots\otimes C$. It follows from the previous proposition that we may form ``polynomials'' in Waldhausen $\infty$-categories (with coefficients in the natural numbers), and the usual formulas hold, such as
\begin{equation*}
( C\otimes D)^m\simeq\bigoplus_{i=0}^m{m\choose i} C^i\otimes D^{m-i}.
\end{equation*}
\end{nul}

\begin{nul} There is a natural generalization \cite[Df. 2.1.1.10]{HA} of the notion of a symmetric monoidal $\infty$-category to that of an \emph{$\infty$-operad}. This is an inner fibration $p:\fromto{O^{\otimes}}{N\Lambda(\FF)}$ satisfying properties that ensure that the objects of $O^{\otimes}$ are, in effect, pairs $(I,X_I)$ consisting of a finite set $I$ and a collection $X_I=\{X_i\}_{i\in I}$ of objects of $O^{\otimes}_{\{\ast\}}$, and that morphisms $\fromto{(J,Y_J)}{(I,X_I)}$ of $C^{\otimes}$ are essentially determined by pairs $(\omega,\phi_I)$ consisting of morphisms $\omega:\fromto{J}{I}$ of $\Lambda(\FF)$ and families of ``multi-morphisms''
\begin{equation*}
\left\{\phi_i:\fromto{Y_{J_i}}{X_i}\right\}_{i\in I}. 
\end{equation*}
An $O^{\otimes}$-algebra in a symmetric monoidal $\infty$-category is simply a morphism of $\infty$-operads \cite[Df. 2.1.2.7]{HA}.
\end{nul}

\begin{dfn} For any $\infty$-operad $O^{\otimes}$, an \textbf{\emph{$O^{\otimes}$-monoidal Waldhausen $\infty$-category}} is an $O^{\otimes}$-algebra in $\Wald_{\infty}^{\otimes}$.
\end{dfn}

\begin{exm} In particular, a \textbf{\emph{monoidal Waldhausen $\infty$-category}} is simply an $O^{\otimes}$-monoidal Waldhausen $\infty$-category, where $O^{\otimes}$ is the associative $\infty$-operad \cite[Df. 4.1.1.3]{HA}. Similarly, a \textbf{\emph{symmetric monoidal Waldhausen $\infty$-category}} will be a $O^{\otimes}$-monoidal Waldhausen $\infty$-category, where $O^{\otimes}$ is the commutative $\infty$-operad \cite[Ex. 2.1.1.18]{HA}.
\end{exm}

\begin{exm} Suppose that $\Phi$ is a perfect operator category \cite{opcat}, and suppose that $\fromto{ X}{N\Lambda(\Phi)}$ is a pair cocartesian fibration such that for any object $I$ of $\Lambda(\Phi)$, the inert morphisms $\fromto{I}{\{i\}}$ induce an equivalence of pairs
\begin{equation*}
\equivto{ X_I}{\prod_{i\in|I|} X_{\{i\}}}.
\end{equation*}
Then we might call $ X$ a \textbf{\emph{$\Phi$-monoidal Waldhausen $\infty$-category}} if, for any object $I$ of $\Lambda(\Phi)$, the pair $ X_I$ is a Waldhausen $\infty$-category, and the morphism
\begin{equation*}
\fromto{\prod_{i\in|I|} X_{\{i\}}\simeq X_I}{ X_{\{\xi\}}}
\end{equation*}
induced by the unique active morphism $\fromto{I}{\{\xi\}}$ is exact separately in each variable. One may show that this notion is essentially equivalent to the notion of $U_{\Phi}^{\otimes}$-monoidal Waldhausen $\infty$-category, where $U_{\Phi}^{\otimes}$ is the symmetrization of the terminal $\infty$-operad over $\Phi$. In particular, an $\OO^{(n)}$-monoidal Waldhausen $\infty$-category is essentially the same thing as an $E_n$-monoidal Waldhausen $\infty$-category.
\end{exm}


\section{Tensor products of virtual Waldhausen $\infty$-categories} The derived $\infty$-category $\mathrm{D}_{\geq0}(k)$ of complexes of vector spaces over a field $k$ with vanishing negative homology inherits a symmetric monoidal structure from the ordinary category of vector spaces. In precisely the same manner, the derived $\infty$-category of Waldhausen $\infty$-categories $\D(\Wald_{\infty})$ --- which is the $\infty$-category of functors $\fromto{\Wald_{\infty}^{\op}}{\SS}$ that preserve all filtered limits and all finite products \cite[Nt. 4.10]{K1} --- inherits a symmetric monoidal structure from $\Wald_{\infty}$.

\begin{prp}\label{prp:VWaldotimes} There exists a symmetric monoidal $\infty$-category $\VWald^{\otimes}$ and a fully faithful symmetric monoidal functor $\into{\Wald_{\infty}^{\otimes}}{\VWald^{\otimes}}$ with the following properties.
\begin{enumerate}[(\ref{prp:VWaldotimes}.1)]
\item The underlying $\infty$-category of $\VWald^{\otimes}$ is the $\infty$-category $\VWald$ of virtual Waldhausen $\infty$-categories, and the underlying functor is the inclusion $\into{\Wald_{\infty}}{\VWald}$.
\item For any symmetric monoidal $\infty$-category $ E^{\otimes}$ whose underlying $\infty$-category admits all sifted colimits, the induced functor
\begin{equation*}
\fromto{\Fun_{N\mathbf{Fin}_{\ast}}(\VWald, E)}{\Fun_{N\mathbf{Fin}_{\ast}}(\Wald_{\infty}, E)}
\end{equation*}
exhibits an equivalence from the full subcategory of spanned by those morphisms of $\infty$-operads $A$ whose underlying functor $A\colon\fromto{\VWald}{ E}$ preserves sifted colimits to the full subcategory of spanned by those morphisms of $\infty$-operads $B$ whose underlying functor $B\colon\fromto{\Wald_{\infty}}{ E}$ preserves filtered colimits.
\item The tensor product functor $\otimes\colon\fromto{\VWald\times\VWald}{\VWald}$ preserves all colimits separately in each variable.
\end{enumerate}
\begin{proof} The only part that is not a consequence of \cite[Pr. 6.3.1.10 and Var. 6.3.1.11]{HA} is the assertion that the tensor product functor
\begin{equation*}
\otimes\colon\fromto{\VWald\times\VWald}{\VWald}
\end{equation*}
preserves direct sums separately in each variable. Pr. \ref{prp:tensorofWaldgoodforsums} states that this holds among Waldhausen $\infty$-categories; the general case follows by exhibiting all the virtual Waldhausen $\infty$-categories as colimits of simplicial diagrams of Waldhausen $\infty$-categories and using the fact that both the tensor product and the direct sum commute with sifted colimits.
\end{proof}
\end{prp}

Since our goal is to study multiplicative structures on additive theories, we should see to it that the symmetric monoidal structure on the derived $\infty$-category of Waldhausen $\infty$-categories gives rise to one on the fissile derived $\infty$-category described in the introduction \ref{nul:add}. It is not the case that the tensor product of two fissile virtual Waldhausen $\infty$-categories is still fissile; however, $\VaddWald$ is the accessible localization of $\VWald$ with respect to the set of morphisms
\begin{equation*}
\left\{i:\fromto{C\oplus C}{E(C)}\ |\ C\in\Wald_{\infty}^{\omega}\right\},
\end{equation*}
where $i$ is the functor defined in \ref{nul:iandr}. We can therefore ask whether the localization functor is \emph{compatible} with the resulting localization functor
\begin{equation*}
L^{\fiss}:\fromto{\VWald}{\VaddWald}.
\end{equation*}
That is, we wish to show that the assignment
\begin{equation*}
\goesto{( X, Y)}{L^{\fiss}( X\otimes Y)};
\end{equation*}
defines a symmetric monoidal structure on $\VaddWald$.

\begin{cnstr} The symmetric monoidal $\infty$-category $\VWald^{\otimes}$ can be described in the following manner. An object $(I, X_I)$ thereof is a finite set $I$ and an $I$-tuple of virtual Waldhausen $\infty$-categories $ X_I\coloneq( X_i)_{i\in I}$. A morphism $\fromto{(J, Y_J)}{(I, X_I)}$ is a morphism $\fromto{J}{I}$ of $\Lambda(\FF)$ along with the data, for every element $i\in I$, of a functor of pairs
\begin{equation*}
\fromto{\bigotimes Y_{J_i}}{ X_i}.
\end{equation*}
We denote by $\VaddWald^{\otimes}$ the full subcategory of $\VWald^{\otimes}$ spanned by those objects $(I, X_I)$ such that for every $i\in I$, $ X_i$ is a distributive virtual Waldhausen $\infty$-category.
\end{cnstr}

\begin{lem}\label{Lfissismultiplicative} The localization functor $L^{\fiss}$ on $\VWald$ of \cite[Pr. 6.7]{K1} is compatible with the symmetric monoidal structure on $\VWald$ is the sense of \cite[Df. 2.2.1.6]{HA}.
\begin{proof} As in \cite[Ex. 2.2.1.7]{HA}, our claim is that for any $L^{\fiss}$-equvalence $\fromto{ U}{ V}$ and any virtual Waldhausen $\infty$-category $ Z$, the morphism $\fromto{ U\otimes Z}{ V\otimes Z}$ is an $L^{\fiss}$-equivalence. Since the tensor product preserves all colimits separately in each variable, we may assume that $ Z$ is a compact Waldhausen $\infty$-category $ D$. Furthermore, since the set of $L^{\fiss}$-equivalences is generated as a strongly saturated class by the set of maps of the form
\begin{equation*}
i\colon\fromto{C\oplus C}{E(C)}
\end{equation*}
in which $C$ is compact, we may assume that $\fromto{U}{ V}$ is of this form. Our claim is thus that for any compact Waldhausen $\infty$-categories $C$ and $D$, the map
\begin{equation*}
i\otimes \id_D:\fromto{(C\oplus C)\otimes D}{E(C)\otimes D}
\end{equation*}
is a $L^{\fiss}$-equivalence. We have a retraction
\begin{equation*}
r\otimes\id_D\colon\fromto{E(C)\otimes D}{(C\oplus C)\otimes D}
\end{equation*}
of this map, and the composition $\fromto{E(C)\otimes D}{E(C)\otimes D}$ is given by the multi-exact functor that carries a pair $(S\ \tikz[baseline]\draw[>=stealth,>->](0,0.5ex)--(0.5,0.5ex);\ T\ \tikz[baseline]\draw[>=stealth,->](0,0.5ex)--(0.5,0.5ex);\ T/S,X)$ to the ``simple tensor''
\begin{equation*}
(S\ \tikz[baseline]\draw[>=stealth,>->](0,0.5ex)--(0.5,0.5ex);\ S\vee(T/S)\ \tikz[baseline]\draw[>=stealth,->](0,0.5ex)--(0.5,0.5ex);\ T/S)\otimes X.
\end{equation*}
That is, the composition $\fromto{E(C)\otimes D}{E(C)\otimes D}$ is given by the direct sum of the functor
\begin{equation*}
\goesto{(S\ \tikz[baseline]\draw[>=stealth,>->](0,0.5ex)--(0.5,0.5ex);\ T\ \tikz[baseline]\draw[>=stealth,->](0,0.5ex)--(0.5,0.5ex);\ T/S,X)}{(S\ \tikz[baseline]\draw[>=stealth,-,double distance=1.5pt](0,0.5ex)--(0.5,0.5ex);\ S\ \tikz[baseline]\draw[>=stealth,->](0,0.5ex)--(0.5,0.5ex);\ 0)\otimes X}
\end{equation*}
and the functor
\begin{equation*}
\goesto{(S\ \tikz[baseline]\draw[>=stealth,>->](0,0.5ex)--(0.5,0.5ex);\ T\ \tikz[baseline]\draw[>=stealth,->](0,0.5ex)--(0.5,0.5ex);\ T/S,X)}{\frac{(S\ \tikz[baseline]\draw[>=stealth,>->](0,0.5ex)--(0.5,0.5ex);\ T\ \tikz[baseline]\draw[>=stealth,->](0,0.5ex)--(0.5,0.5ex);\ T/S)\otimes X}{(S\ \tikz[baseline]\draw[>=stealth,-,double distance=1.5pt](0,0.5ex)--(0.5,0.5ex);\ S\ \tikz[baseline]\draw[>=stealth,->](0,0.5ex)--(0.5,0.5ex);\ 0)\otimes X}}.
\end{equation*}
But in $\VaddWald$, this is homotopic to the identity.
\end{proof}
\end{lem}

\begin{nul} Note that as a corollary, we find that, for any integer $m\geq 0$, we obtain $L^{\fiss}$-equivalences
\begin{equation*}
 F_m( C)\otimes D\simeq m C\otimes D\simeq F_m( C\otimes D),
\end{equation*}
where $F_m(C)$ is the Waldhausen $\infty$-category of filtered objects $X_0\ \tikz[baseline]\draw[>=stealth,>->](0,0.5ex)--(0.5,0.5ex);\ \cdots\ \tikz[baseline]\draw[>=stealth,>->](0,0.5ex)--(0.5,0.5ex);\ X_m$ \cite[Nt. 5.5]{K1}. Consequently, we obtain an $L^{\fiss}$-equivalence
\begin{equation*}
 F_m( C)\simeq F_m(N\mathbf{Fin}_{\ast})\otimes C.
\end{equation*}

This observation yields another way to think about the suspension functor in $\VaddWald$. Just as suspension is smashing with a circle in the homotopy theory of spaces, we have
\begin{equation*}
\Sigma C\simeq S( C)\simeq S(N\mathbf{Fin}_{\ast}\otimes C)\simeq S(N\mathbf{Fin}_{\ast})\otimes L^{\fiss}( C)
\end{equation*}
in $\VaddWald$. Here the fissile virtual Waldhausen $\infty$-category $ S(N\mathbf{Fin}_{\ast})$ is playing the role of a circle. The algebraic $K$-theory of a Waldhausen $\infty$-category $ C$ can thus be described as the space
\begin{equation*}
\Omega I S( C)\simeq\Omega I S(N\mathbf{Fin}_{\ast}\otimes C)\simeq\Omega I( S(N\mathbf{Fin}_{\ast})\otimes L^{\fiss}( C)),
\end{equation*}
where $I$ is the left derived functor of $\iota$. More generally, for any pre-additive theory $\phi$, the additivization $D\phi$ can be computed by the formula
\begin{equation*}
D\phi( C)\simeq\Omega\Phi S( C)\simeq\Omega\Phi S(N\mathbf{Fin}_{\ast}\otimes C)\simeq\Omega\Phi( S(N\mathbf{Fin}_{\ast})\otimes L^{\fiss}( C)),
\end{equation*}
where $\Phi$ is the left derived functor of $\phi$.
\end{nul}

\begin{prp}\label{prp:Dfissotimes} The functor $\fromto{\VaddWald^{\otimes}}{N\Lambda(\FF)}$ is a symmetric monoidal $\infty$-category with the property that the tensor product
\begin{equation*}
\otimes\colon\fromto{\VaddWald\times\VaddWald}{\VaddWald}
\end{equation*}
preserves all colimits separately in each variable.
\begin{proof} That $\VaddWald^{\otimes}$ is symmetric monoidal follows from \cite[Pr. 2.2.1.9]{HA}. Observe that the functor
\begin{equation*}
\otimes\colon\fromto{\VaddWald\times\VaddWald}{\VaddWald}
\end{equation*}
can be identified with the functor given by the assignment
\begin{equation*}
\goesto{( X, Y)}{L^{\fiss}( X\otimes Y)};
\end{equation*}
hence it follows from \cite[Cor. 6.7.2]{K1} it preserves colimits in each variable.
\end{proof}
\end{prp}

\begin{cnstr} The stabilization
\begin{equation*}
\Sp(\VaddWald)
\end{equation*}
of the fissile derived $\infty$-category $\VaddWald$ inherits a canonical symmetric monoidal structure $\Sp(\VaddWald)^{\otimes}$, given by applying the symmetric monoidal functor
\begin{equation*}
L^{\otimes}\colon\fromto{{\Pr}^{\mathrm{L},\otimes}}{{\Pr}^{\mathrm{L},\otimes}_{\mathrm{St}}}
\end{equation*}
induced by the localization $L=\Sp\otimes-$ of \cite[Pr. 6.3.2.17]{HA} to the presentable $\infty$-category $\VaddWald$. By construction, the tensor product functor $\otimes$ preserves colimits separately in each variable. Furthermore, the functor
\begin{equation*}
\Sigma^{\infty}\colon\fromto{\VaddWald}{\Sp(\VaddWald)}
\end{equation*}
is symmetric monoidal.

We call the $\infty$-category $\Sp(\VaddWald)$ the \textbf{\emph{stable $\infty$-category of Waldhausen $\infty$-categories}}, and we call its homotopy category $h\Sp(\VaddWald)$ the \textbf{\emph{stable homotopy category of Waldhausen $\infty$-categories}}.

Note that $\Sp(\VaddWald)$ is equivalent to the full subcategory of the $\infty$-category $\Fun(\Wald_{\infty}^{\omega,\op},\Sp)$ spanned by those functors $ X\colon\fromto{\Wald_{\infty}^{\omega,\op}}{\Sp}$ such that for any compact Waldhausen $\infty$-category $C$, the morphisms induced by the exact functor $i\colon\fromto{C\oplus C}{E(C)}$ exhibits $X(E(C))$ as the direct sum of $X(C)$ and $X(C)$.
\end{cnstr}

\begin{nul} Observe that the $K$-groups of a Waldhausen $\infty$-category $C$ are given by $\Ext$ groups in the stable homotopy category of Waldhausen $\infty$-categories. Indeed, just as one may describe the stable homotopy groups of a pointed space $X$ as Ext groups out of the sphere spectrum:
\begin{equation*}
\pi_n^{s}(X)\cong\Ext^{-n}_{\Sp}(\Sigma^{\infty}S^0,\Sigma^{\infty}X),
\end{equation*}
so too may one describe the algebraic $K$-theory groups of a Waldhausen $\infty$-category $C$ as Ext groups out of the suspension spectrum of $\Fin_{\ast}$:
\begin{equation*}
K_n(C)\cong\Ext^{-n}_{\Sp(\VaddWald)}(\Sigma^{\infty}\Fin_{\ast},\Sigma^{\infty}C).
\end{equation*}
\end{nul}


\section{Multiplicative theories} Now we are in a position to study theories that are compatible with the monoidal structure on Waldhausen $\infty$-categories.

\begin{nul} Recall \cite[Df. 7.1]{K1} that for any $\infty$-topos $ E$, an $ E$-valued theory is a reduced functor
\begin{equation*}
\phi\colon\fromto{\Wald_{\infty}}{ E}
\end{equation*}
that preserves filtered colimits.
\end{nul}

\begin{dfn} Suppose $ E$ an $\infty$-topos. A \textbf{\emph{multiplicative theory valued in $ E$}} is a morphism of $\infty$-operads
\begin{equation*}
\phi^{\otimes}\colon\fromto{\Wald_{\infty}^{\otimes}}{ E^{\times}}
\end{equation*}
such that the underlying functor $\phi\colon\fromto{\Wald_{\infty}}{ E}$ preserves all filtered colimits and carries the zero object to the terminal object. We will say that the multiplicative theory $\phi^{\otimes}$ \textbf{\emph{extends}} the theory $\phi$.

We denote by
\begin{equation*}
\Thy^{\otimes}( E)\subset\Fun_{N\Lambda(\FF)}(\Wald_{\infty}^{\otimes}, E^{\times})
\end{equation*}
the full subcategory spanned by the multiplicative theories, and we denote by $\Add^{\otimes}( E)$ the full subcategory spanned by those multiplicative theories such that the underlying theory $\fromto{\Wald_{\infty}}{ E}$ is additive.
\end{dfn}

\begin{nul} It follows from Pr. \ref{prp:VWaldotimes} that any multiplicative theory
\begin{equation*}
\phi^{\otimes}\colon\fromto{\Wald_{\infty}^{\otimes}}{ E^{\times}}
\end{equation*}
can be extended to a reduced functor $\Phi^{\otimes}\colon\fromto{\VWald^{\otimes}}{ E^{\times}}$ of $\infty$-operads such that the underlying functor $\fromto{\VWald}{ E_{\ast}}$ preserves sifted colimits. This is the \textbf{\emph{multiplicative left derived functor}} of $\phi^{\otimes}$.
\end{nul}

\begin{exm} The theory $\iota\colon\fromto{\Wald_{\infty}}{\Kan}$ \cite[Nt. 1.7]{K1} can be extended to a multiplicative theory $\iota^{\otimes}$ in the following manner. The right adjoint $\fromto{\Cat_{\infty}}{\Kan}$ of the inclusion can be given the structure of a symmetric monoidal functor
\begin{equation*}
\fromto{\Cat_{\infty}^{\times}}{\Kan^{\times}}
\end{equation*}
for the cartesian symmetric monoidal structures in an essentially unique manner. The desired morphism $\iota^{\otimes}\colon\fromto{\Wald_{\infty}^{\otimes}}{\Kan^{\times}}$ of $\infty$-operads is now the composite
\begin{equation*}
\Wald_{\infty}^{\otimes}\ \tikz[baseline]\draw[>=stealth,right hook->](0,0.5ex)--(0.5,0.5ex);\ \Pair_{\infty}^{\times}\to\Cat_{\infty}^{\times}\to\Kan^{\times}.
\end{equation*}
\end{exm}

We now wish to show that additivizations of multiplicative theories are naturally multiplicative (Th. \ref{thm:main}). Since additivizations are constructed via Goodwillie differentials, we will prove a general result about these. What follows is surely not the most general result one can prove, but it's enough for our purposes.

\begin{prp}\label{prp:diffofalg} Suppose $C^{\otimes}$ and $D^{\otimes}$ symmetric monoidal $\infty$-categories. Suppose that the underlying $\infty$-category $D$ is presentable, and suppose that the underlying $\infty$-category $C$ is small and that it admits a terminal object and all finite colimits. Finally, assume that $C^{\otimes}$ is compatible with all finite colimits, and assume that $D^{\otimes}$ is compatible with all colimits. Then the inclusion
\begin{equation*}
\into{\Exc(C,D)\times_{\Fun(C,D)}\Alg_{C^{\otimes}}(D^{\otimes})}{\Alg_{C^{\otimes}}(D^{\otimes})}
\end{equation*}
admits a left adjoint.
\begin{proof}[Proof A] The Adjoint Functor Theorem \cite[Cor. 5.5.2.9]{HTT} and the existence of excisive approximation \cite[Th. 7.1.10]{HA} implies that $\Exc(C,D)\subset\Fun(C,D)$ is stable under both limits and $\kappa$-filtered colimits for some regular cardinal $\kappa$. The fiber product
\begin{equation*}
\Exc(C,D)\times_{\Fun(C,D)}\Alg_{C^{\otimes}}(D^{\otimes})
\end{equation*}
can be identified with the full subcategory of $\Alg_{C^{\otimes}}(D^{\otimes})$ spanned by those functors $F^{\otimes}\colon\fromto{C^{\otimes}}{D^{\otimes}}$ over $N\Lambda(\FF)$ with the property that the underlying functor $F\colon\fromto{C}{D}$ is excisive. It now follows from \cite[Cor. 3.2.2.5 and Pr. 3.2.3.1(4)]{HA} that this subcategory is stable under both limits and $\kappa$-filtered colimits. A second coat of the Adjoint Functor Theorem \cite[Cor. 5.5.2.9]{HTT} now yields the result.
\end{proof}

\begin{proof}[Proof B] Alternately, one may prove this result (in fact a more general version thereof) in an explicit fashion by applying a variant of \cite[Cor. 5.2.7.11]{HTT}. Note that for the argument there to go through, one does not need the full strength of the condition that the forgetful functor
\begin{equation*}
p\colon\fromto{\Alg_{C^{\otimes}}(D^{\otimes})}{\Fun(C,D)}
\end{equation*}
be a cocartesian fibration. (It obviously isn't in our case.) One need only ensure that for every multiplicative functor $F^{\otimes}\colon\fromto{C^{\otimes}}{D^{\otimes}}$ with underlying functor $F\colon\fromto{C}{D}$, there exist a localization $\theta\colon\fromto{F}{P_1F}$ with respect to $\Exc(C,D)$ and a $p$-cocartesian edge $\fromto{F^{\otimes}}{(P_1F)^{\otimes}}$ that covers $\theta$. Recall \cite[Cnstr. 7.1.1.27]{HA} that $P_1F$ can be obtained as the sequential colimit of the sequence
\begin{equation*}
F\to T_1F\to T_1^2F\to \cdots,
\end{equation*}
where $T_1F=\Omega\circ F\circ\Sigma$ as in \ref{nul:aboutP1}.

Let $\NN^{+}$ denote the following ordinary category. An object is a pair $(I,k_I)$ consisting of a finite set $I$ and an $I$-tuple $k_I=(k_i)_{i\in I}$ of natural numbers. A morphism $\fromto{(J,\ell_J)}{(I,k_I)}$ is a map of finite sets $\fromto{J}{I_{+}}$ such that for any element $i\in I$, one has
\begin{equation*}
\sum_{j\in J_i}\ell_j\leq k_i.
\end{equation*}
The forgetful functor $\fromto{\NN^{+}}{\Lambda(\FF)}$ is the Grothendieck opfibration that corresponds to the functor $\fromto{\Lambda(\FF)}{\Cat}$ that exhibits the ordered set of natural numbers as a symmetric monoidal category under addition. Hence the cocartesian fibration $\fromto{N\NN^{+}}{N\Lambda(\FF)}$ is a symmetric monoidal $\infty$-category. Now one may define a functor
\begin{equation*}
(P_1F)^{\boxtimes}\colon\fromto{N\NN^{+}\times_{N\mathbf{Fin}_{\ast}}C^{\otimes}}{D^{\otimes}}
\end{equation*}
such that the formula
\begin{equation*}
(P_1F)^{\boxtimes}(I,k_I,X_I)\coloneq(T_1^{k_i}X_i)_{i\in I}
\end{equation*}
holds. Then our desired functor $(P_1F)^{\otimes}$ will be the left Kan extension of $(P_1F)^{\boxtimes}$ along the projection $\fromto{N\NN^{+}\times_{N\Lambda(\FF)}C^{\otimes}}{C^{\otimes}}$.

We also have the functor $F^{\boxtimes}\colon\fromto{N\NN^{+}\times_{N\Lambda(\FF)}C^{\otimes}}{D^{\otimes}}$, which is the projection $\fromto{N\NN^{+}\times_{N\Lambda(\FF)}C^{\otimes}}{C^{\otimes}}$ composed with $F^{\otimes}\colon\fromto{C^{\otimes}}{D^{\otimes}}$. The natural transformations $\fromto{F}{T_1^kF}$ induce a natural transformation $\fromto{F^{\boxtimes}}{(P_1F)^{\boxtimes}}$ and thus an induced natural transformation $\theta^{\otimes}\colon\fromto{F^{\otimes}}{(P_1F)^{\otimes}}$. A quick computation shows that $\theta^{\otimes}$ is a cocartesian edge covering $\theta$.
\end{proof}
\end{prp}

\begin{cor}\label{cor:diffofalg} Suppose $C^{\otimes}$ and $D^{\otimes}$ symmetric monoidal $\infty$-categories. Suppose that the underlying $\infty$-category $D$ is presentable, and suppose that the underlying $\infty$-category $C$ is compactly generated. Finally, assume that $C^{\otimes}$ and $D^{\otimes}$ are compatible with all colimits. Then the inclusion
\begin{equation*}
\into{\Exc_{ F}(C,D)\times_{\Fun(C,D)}\Alg_{C^{\otimes}}(D^{\otimes})}{\Alg_{C^{\otimes}}(D^{\otimes})}
\end{equation*}
admits a left adjoint, where $\Exc_{ F}(C,D)$ denotes the full subcategory of $\Fun(C,D)$ spanned by excisive functors $\fromto{C}{D}$ that preserve all filtered colimits.
\end{cor}

With this result in hand, we easily prove our main theorem.

\begin{thm}\label{thm:main} The inclusion $\into{\Add^{\otimes}( E)}{\Thy^{\otimes}( E)}$ admits a left adjoint that covers the left adjoint of the inclusion $\into{\Add( E)}{\Thy( E)}$.
\begin{proof} In light of \cite[Th. 7.6]{K1}, the claim is that the inclusion
\begin{equation*}
\begin{tikzpicture} 
\matrix(m)[matrix of math nodes, 
row sep=3ex, column sep=3ex, 
text height=1.5ex, text depth=0.25ex] 
{\Exc_{\ast}^{ G}(\VaddWald, E)\times_{\Fun(\VaddWald, E)}\Alg_{\VaddWald^{\otimes}}( E^{\times})\\ 
\Fun_{\ast}^{ G}(\VaddWald, E)\times_{\Fun(\VaddWald, E)}\Alg_{\VaddWald^{\otimes}}( E^{\times})\\}; 
\path[>=stealth,right hook->,font=\scriptsize] 
(m-1-1) edge (m-2-1); 
\end{tikzpicture}
\end{equation*}
admits a left adjoint. We now appeal to Cor. \ref{cor:diffofalg}, and the proof is completed by the observation that if a functor $\fromto{\VaddWald}{ E}$ preserves geometric realizations, then so does Goodwillie differential. (This follows from \cite[Lm. 7.7]{K1}; see the proof of \cite[Th. 7.8]{K1}.)
\end{proof}
\end{thm}

\begin{ntn} Write $D^{\otimes}$ for the left adjoint of the inclusion
\begin{equation*}
\into{\Add^{\otimes}( E)}{\Thy^{\otimes}( E)}.
\end{equation*}
\end{ntn}

We may now define $K^{\otimes}\coloneq D^{\otimes}\iota^{\otimes}$, whence we deduce that algebraic $K$-theory is naturally a multiplicative theory.
\begin{prp} There exists a canonical multiplicative extension $K^{\otimes}$ of algebraic $K$-theory.
\end{prp}

In light of \cite[Pr. 7.2.4.14 and Pr. 7.2.6.2]{HA}, we also find the following.

\begin{cor} There exists a canonical multiplicative extension $\KK^{\otimes}$ of the connective algebraic $K$-theory functor $\KK\colon\fromto{\Wald_{\infty}}{\Sp}$.
\end{cor}

\begin{cor}\label{cor:Kpreservesalg} For any $\infty$-operad $O^{\otimes}$, composition with the multiplicative extensions $K^{\otimes}$ and $\KK^{\otimes}$ induce functors
\begin{equation*}
K^{\otimes}\colon\fromto{\Alg_{O^{\otimes}}(\Wald_{\infty}^{\otimes})}{\Alg_{O^{\otimes}}(\Kan^{\times})}
\end{equation*}
and
\begin{equation*}
\KK^{\otimes}\colon\fromto{\Alg_{O^{\otimes}}(\Wald_{\infty}^{\otimes})}{\Alg_{O^{\otimes}}(\Sp^{\wedge})}.
\end{equation*}
\end{cor}
\noindent As a special case of this, we obtain the following.
\begin{cor} Suppose $ A^{\otimes}$ a monoidal Waldhausen $\infty$-category. Then composition with the multiplicative extensions $K^{\otimes}$ and $\KK^{\otimes}$ induce functors
\begin{equation*}
K^{\otimes}\colon\fromto{\mathrm{LMod}_{ A^{\otimes}}(\Wald_{\infty}^{\otimes})}{\mathrm{LMod}_{K^{\otimes}( A^{\otimes})}(\Kan^{\times})}
\end{equation*}
and
\begin{equation*}
\KK^{\otimes}\colon\fromto{\mathrm{LMod}_{ A^{\otimes}}(\Wald_{\infty}^{\otimes})}{\mathrm{LMod}_{\KK^{\otimes}( A^{\otimes})}(\Sp^{\wedge})}.
\end{equation*}
\end{cor}

\begin{exm}[Deligne Conjecture for algebraic $K$-theory]\label{ex:Deligne} Suppose $C^{\otimes}$ a pointed, symmetric monoidal $\infty$-category that admits all finite colimits. Assume that the tensor product $\otimes\colon\fromto{C\times C}{C}$ preserves finite colimits separately in each variable. Then $C^{\otimes}$ may be viewed as an $E_{\infty}$ object of $\Wald_{\infty}^{\otimes}$, and the $K$-theory spectrum of $C$ is naturally endowed with an $E_{\infty}$-structure given by $\KK^{\otimes}(C^{\otimes})$ [Cor. \ref{cor:Kpreservesalg}].

Carrying an $E_n$-algebra in $C^{\otimes}$ to its $\infty$-category of right modules defines a functor $\Theta_{ C}$ \cite[Rk. 6.3.5.15]{HA}, which factors through a functor
\begin{equation*}
\Alg_{E_n}(C^{\otimes})\to\Alg_{E_{n-1}}(\Wald_{\infty}^{\otimes}).
\end{equation*}
Composing this with $K^{\otimes}$ and $\KK^{\otimes}$, we obtain functors
\begin{equation*}
\fromto{\Alg_{E_n}(C^{\otimes})}{\Alg_{E_{n-1}}(\Kan^{\times})}\textrm{\quad and\quad}\fromto{\Alg_{E_n}(C^{\otimes})}{\Alg_{E_{n-1}}(\Sp^{\wedge})}.
\end{equation*}
\end{exm}

\begin{exm} The previous example makes an assortment of \emph{iterated $K$-theories} possible. If $n\geq 1$, then by forming iterated compositions of the various functors $\fromto{\Alg_{E_k}(\Sp^{\wedge})}{\Alg_{E_{k-1}}(\Sp^{\wedge})}$ constructed above, we obtain \textbf{\emph{$n$-fold algebraic $K$-theory functor}}
\begin{equation*}
\KK^{(n)}\colon\fromto{\Alg_{E_n}(\Sp^{\wedge})}{\Sp^{\wedge}}
\end{equation*}
as well as an infinite hierarchy of functors
\begin{equation*}
\KK^{(n)}\colon\fromto{\Alg_{E_\infty}(\Sp^{\wedge})}{\Alg_{E_\infty}(\Sp^{\wedge})}.
\end{equation*}
The \emph{Chromatic Red Shift Conjectures} of Ausoni and Rognes (presaged by Waldhausen and Hopkins) implies that $\KK^{(n)}$ should, in effect, carry $E_n$-rings of telescopic complexity $m$ to spectra of telescopic complexity $m+n$. We hope to investigate these phenomena in future work.
\end{exm}

\begin{exm} We also find that Waldhausen's $A$-theory of an $n$-fold loopspace $X$ (defined as in \cite[Ex. 2.10]{K1}) carries a canonical $E_{n-1}$-monoidal structure.
\end{exm}

\begin{nul} Note that, although this result ensures that algebraic $K$-theory is merely \emph{lax} symmetric monoidal, as a functor to spectra it's actually slightly better: by Barratt--Priddy--Quillen, $\KK^{\otimes}\colon\fromto{\Wald_{\infty}^{\otimes}}{\Sp^{\wedge}}$ carries the unit $N\mathbf{Fin}_{\ast}$ to the unit $S^{0}$. 
\end{nul}

Let us conclude by appealing to the recent work \cite{GlasmanDay} of Saul Glasman. Glasman identifies the $\infty$-category
$\Alg_{\Wald^{\otimes}}( E^{\times})$ with the $\infty$-category of $E_{\infty}$ algebras in $\Fun(\Wald_{\infty}, E)$, equipped with the Day convolution structure. When $ E=\Kan$, one sees that the functor $\iota$ is corepresented by the unit $N\mathbf{Fin}_{\ast}$; consequently it is the unit for the Day convolution product. Hence it admits a unique $E_{\infty}$ structure, which under Glasman's equivalence must coincide with the multiplicative theory $\iota^{\otimes}$. Moreover, $\iota^{\otimes}$ is initial in $\Alg_{\Wald^{\otimes}}( E^{\times})$. Meanwhile, the universal property of $K^{\otimes}$ ensures that for any additive multiplicative theory $\phi^{\otimes}$, we have
\begin{equation*}
\Map(K^{\otimes},\phi^{\otimes})\simeq\Map(\iota^{\otimes},\phi^{\otimes})\simeq\ast.
\end{equation*}
That is, there is an essentially unique multiplicative ``trace map'' from $K$-theory to \emph{any} additive and multiplicative theory. In other words, we have the following.

\begin{thm}\label{multtrace} Algebraic $K$-theory of Waldhausen $\infty$-categories is the initial additive and multiplicative theory.
\end{thm}


\bibliographystyle{amsplain}
\bibliography{kthy}

\end{document}